\documentstyle[epsf,amssymb,11pt]{article}
\newtheorem{theorem}{Theorem}

\newtheorem{lemma}[theorem]{Lemma}

\newtheorem{conjecture}[theorem]{Conjecture}
\newtheorem{corollary}[theorem]{Corollary}

\newcommand{\ie}{{\em i.e.}}
\newcommand{\Tr}{\mbox{Tr}}
\newcommand{\tensor}{\otimes}

\newcommand{\C}{{\Bbb C}}
\newcommand{\Z}{{\Bbb Z}}
\newcommand{\Sz}{\mbox{Sz}}
\newcommand{\Aut}{\mbox{Aut}}

\newcommand{\fig}[2]{
    \begin{figure}[ht]
    \centerline{\epsfbox{#1.eps}}
    \caption{#2}\label{#1}
    \end{figure}}

\newcommand{\lc}{\raisebox{-.09in}{\makebox[.5in]{\epsfbox{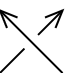}}}}
\newcommand{\rc}{\raisebox{-.09in}{\makebox[.5in]{\epsfbox{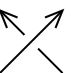}}}}
\newcommand{\mc}{\raisebox{-.09in}{\makebox[.5in]{\epsfbox{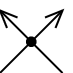}}}}

\begin{document}

\title{Detecting knot invertibility}
\author{Greg Kuperberg}
\date{February 21, 1995}
\maketitle
\begin{abstract}
We discuss the consequences of the possibility that Vassiliev invariants do
not detect knot invertibility as well as the fact that quantum Lie group
invariants are known not to do so.  On the other hand, finite group
invariants, such as the set of homomorphisms from the knot group to $M_{11}$,
can detect knot invertibility.  For many natural classes of knot invariants,
including Vassiliev invariants and quantum Lie group invariants, we
can conclude that the invariants either distinguish all oriented knots, or
there exist prime, unoriented knots which they do not distinguish.
\end{abstract}

Not long ago, two optimistic conjectures about the nature
of Vassiliev invariants were widely circulated:

\begin{conjecture} Vassiliev invariants distinguish all prime, unoriented
knots.~\label{cvtok}
\end{conjecture}

\begin{conjecture} \cite{Bar-Natan:vassiliev} All Vassiliev invariants are
linear combinations of derivatives at $q = 1$ of quantum Lie group
invariants.~\label{cqtov}
\end{conjecture}

Let $s$ be a satellite operation on knots, so that if $K$ is a knot,
$s(K)$ is a satellite of $K$.  We say that an $X$-valued invariant $V$ on knots
{\em intertwines} $s$ if $V(s(K)) = s_X(V(K))$ for some $s_X:X \to X$, so
that the invariant $V \circ s$ yields no more information than $V$
itself.  In this paper, we present the following three
theorems.

\begin{theorem} Let $V$ be a knot invariant which intertwines satellite
operations.  Then either $V$ distinguishes all oriented knots, or there exist
prime, unoriented knots which are not distinguished by $V$.\label{thsat}
\end{theorem}

More precisely, Theorem~\ref{thsat} states that, assuming $V$ distinguishes
prime, unoriented knots, it distinguishes any pair $K_1$ and $K_2$ of
distinct oriented knots.  Strictly speaking, we will only prove this when
$K_2$ is the inverse of $K_1$ (\ie, $K_1$ with the opposite orientation),
but this is in fact the main case.  The alternative case where they are not
inverses is settled as follows:  Knot inversion is itself a satellite
operation, so $V(K_1) = V(K_2)$ implies $V(-K_1) = V(-K_2)$, which is to say
that $V$ does not distinguish $K_1$ and $K_2$ as unoriented knots. If they
are not both prime, then by Lemma~\ref{lprime}, they can be made prime by a
satellite operation and they will remain indistinguishable.

\begin{theorem} The universal Vassiliev invariant $v_n$ of order $n$
intertwines satellite operations. \label{thvass}
\end{theorem}

\begin{theorem} The universal quantum link invariant far a Lie
algebra $\frak g$ intertwines satellite operations. \label{thquantum}
\end{theorem}

It is known that the quantum link invariants do not detect invertibility.
Therefore, the three theorems have the following corollary.

\begin{corollary} Conjectures~\ref{cvtok} and~\ref{cqtov} are mutually
exclusive.
\end{corollary}

Our arguments are inspired by the close relation between inverting a knot
and the operation of mild mutation on a 3-manifold $M$ such as a knot
complement.  Mild mutation consists of cutting $M$ along a torus
and then regluing after applying the central involution of the
mapping class group.  This operation was similarly exploited by
Kania-Bartoszynska to exhibit two closed 3-manifolds with the same
quantum 3-manifold invariants \cite{Kania-Bartoszynska}.

Despite the pessimistic aspect of Theorem~\ref{thsat}, it seems likely
that quantum invariants, which for links are determined by Vassiliev invariants,
distinguish all atoroidal, irreducible 3-manifolds, in particular unoriented knots
which are not satellites.

\section{Satellite operations}

Let $\cal K$ be the set of framed, oriented knots in $S^3$ up to isotopy.
Suppose that $L$ is a connected, framed, oriented two-component link whose
components are labelled as $L_1$ and $L_2$.  Suppose further that $L_1$ is an
un-knot with untwisted framing. Then a knot $K \in \cal K$ together with $L$
yield a satellite knot $s_L(K)$ of $K$ by removing a tube around $K$ and a
tube around $L_1$ and gluing the two knot complements together to retain only
$L_2$, specifically by gluing the longitude of $K$ (given by the framing) to
the meridian of $L_1$ and vice-versa.  In other words, $L$ yields a satellite
operation $s_L:{\cal K} \to {\cal K}$.

There are two ways in which $s_L(K)$ might be a composite knot.
Firstly, any summand of $L$ (necessarily a summand of $L_2$) appears as
a summand of $s_L(K)$.  Secondly, if $K$ is trivial, then
$s_L(K)$ can be substantially simpler than and different from $L$
and can be a composite knot, even if $L$ is prime.

\begin{lemma} If $L$ is prime and $K$ is non-trivial, then 
$s_L(K)$ is prime, even if $K$ is composite. \label{lprime}
\end{lemma}

This is a classical result that is proved by an innermost circle
argument:  Assume a sphere $S$ that separates $s_L(K)$ into summands
and put it in a position of minimal transverse intersection with
the boundary torus of $K$.  An analysis of the innermost circles
of intersection on $S$ yields a summand of $L$ or a new sphere
with smaller intersection.

Now suppose that $L$ is hyperbolic and does not admit a symmetry that inverts
$L_1$. Let $K$ be a non-invertible knot and let $-K$ be its inverse. Then
$s_L(K)$ and $s_L(-K)$ are prime, and they differ as unoriented knots as
follows:  Since $L$ is hyperbolic, its complement is an aspherical,
atoroidal, and acylindrical 3-manifold (\ie, $L$ is prime, connected, not the
Hopf link, and not a non-trivial satellite).  Let $T$ be the torus jacketing
$s_L(K)$ (\ie, the gluing torus in the construction of the satellite), and
let $h$ be a diffeomorphism of $S^3$ that takes $s_L(K)$ to $s_L(-K)$ such
that $T$ and $h(T)$ are in a position of minimal transverse intersection.  If
the intersection is non-empty, then by another standard analysis of innermost
circles, the result is either an essential sphere, cylinder, or torus, or a
smaller intersection \cite{Jaco-Shalen}.  If the intersection is empty, then
$T$ and $h(T)$ must be parallel, for otherwise $L$ would have an essential
torus.  But in this case, $h$ must invert either the $K$ piece or the $L$
piece, a contradiction.

For an explicit example, we can take $K$ to be the knot $8_{17}$, as shown in
Figure~\ref{f817}, and $L$ to be the link $9^2_{34}$, as shown in
Figure~\ref{f9234}, using Conway's enumeration \cite{Rolfsen}.  According to
SnapPea \cite{SnapPea} and Mostow rigidity, neither $K$ nor $L$ admit the
undesired symmetries.  Indeed, $L$ has no symmetries, and we can take either
component to be $L_1$.

\fig{}{The non-invertible knot $8_{17}$.}
\fig{}{The link $9^2_{34}$, which has no symmetries.}

Suppose that $X$ is a set and $f:{\cal K} \to X$ is an invariant which
intertwines all satellite operations.  We say that $f$ {\em distinguishes}
$K_1$ and $K_2$ as unoriented knots if the sets $\{f(K_1),f(-K_1)\}$ and
$\{f(K_2),f(-K_2)\}$ differ.  Let $K$ and $L$ be as above, let $L'$ be $L$
with $L_2$ inverted, and suppose that $f$ distinguishes $s_L(K)$ and
$s_L(-K)$ as unoriented knots. Then it is immediate that either $f \circ s_L$
or $f \circ s_{L'}$ distinguishes $K$ from $-K$.  Since $f$ intertwines
these, $f$ must also distinguish $K$ from $-K$. This establishes
Theorem~\ref{thsat}.

It remains to show that $v_n$, the universal Vassiliev invariant of order $n$,
and $Q_{q,{\frak g}}(K)$, the universal quantum invariant for a knot $K$ and
a Lie algebra $\frak g$, both intertwine all satellite operations.

\section{Vassiliev invariants}

As before, let ${\cal K}_0 = {\cal K}$ be the set of framed, oriented knots in
$S^3$, and let ${\cal K}_n$ be the set of framed, oriented, immersed circles
in $S^3$ with $n$ double points with transverse tangents.  (The immersed
circles will loosely be called generalized knots.)  For example, $\cal K$ and
${\cal K}_1$ are the codimension 0 and codimension 1 cells of a stratification
of the space of smooth maps from $S^1$ to $S^3$.  Let $V$ be the free abelian
group generated by the union of all ${\cal K}_n$'s, modulo the relation that
$$\mc = \rc - \lc$$
for every triple of knots which differ at one crossing as indicated. Note that
$\cal K$ is a basis for $V$.  The Vassiliev space $V_n$ of order $n$ is the
quotient of $V$ by the subgroup generated by ${\cal K}_{n+1}$, and the
universal Vassiliev invariant of order $n$ is the induced map $v_n:{\cal K}
\to V_n$.  Theorem~\ref{thvass} asserts that for any two-component link $L$,
$v_n \circ s_L = s_{L,V_n} \circ v_n$ for some endomorphism $s_{L,V_n}$ of
$V_n$.  Shifting $n$ by 1, if $J_n$ is the subspace (or ideal) of $V$
generated by ${\cal K}_n$, then, equivalently, we wish to show
that $s_L(J_n) \subseteq J_n$.

Let $C$ be the free abelian group whose basis consists of a left-handed
crossing $l$ and a right-handed crossing $r$, not part of any knot:
$$l = \lc$$
$$r = \rc$$
Let $a = r - l \in C$, and let $A_{1,1}$ be the subgroup of $C$ spanned by
$a$. More generally, let $A_{n,k} \subset C^{\tensor n}$ be the span of all
subgroups of the form $A_{1,1}^{\tensor k} \tensor C^{\tensor n-k}$ with the
tensor factors permuted arbitrarily.  Given a generalized knot $K \in {\cal
K}_n$ with linearly ordered double points (not necessarily
ordered as they appear along $K$), there is a map $\phi_K:C^{\tensor
n} \to V$ which acts on an element of the product basis by replacing
each double point of $K$ by the corresponding factor crossing.  For example,
if $K$ is the generalized knot shown in Figure~\ref{fgenknot}, then $\phi_K(l
\tensor r)$ is the knot shown in Figure~\ref{fphi}. In
particular, $\phi_K(a^{\tensor n})$ is simply $K$ itself as an element of
$V$. More generally, $\phi_K(A_{n,k}) \subseteq {\cal K}_k$. Finally, let
$a_{n,\Delta} = r^{\tensor n} - l^{\tensor n} \in C^{\tensor n}$.

\fig{}{A generalized knot $K$ with ordered double points.}
\fig{}{The knot $\phi_K(l \tensor r)$.}

Let $L$ be a 2-component link and let $K$ be as before.  We can take the
satellite of $K$, of sorts, corresponding to $L$ by replacing every double
point by $m^2$ double points, replacing every crossing by $m^2$ crossings,
every arc but one by $m$ arcs, and the remaining arc by the pattern of $L_2$,
as shown in Figure~\ref{fl2}. The result $K'$ is not unique; it depends on
where we place the pattern of $L_2$, and if $L$ is rearranged, even $n$ can
change. Nevertheless, the vector
$$b = \phi_{K'}(a_{m^2,\Delta}^{\tensor n})$$
in $V$ is uniquely determined by $K$ and $L$ up to sign, because it equals
$\pm s_L(K)$, where $s_L$ is interpreted as an endomorphism of $V$ and $K$ is
interpreted as an element of $V$. Here we order the double points of
$K'$ in bunches according to which double points of $K$ they come from. (The
sign is due to the fact that the satellite operations may replace a
left-handed crossing by a combination of right- and left-handed crossings.)
In other words, $K$ represents a certain linear combination of honest knots
in $V$, and $s_L(K) = \pm b$ is the same linear combination of the
$L$-satellites of these knots.

\fig{}{A fake satellite operation on generalized knots.}

Observe that $a_{m^2,\Delta} \in A_{m^2,1}$.  It follows that 
$$a_{m^2,\Delta}^{\tensor n} \in A_{m^2,1}^{\tensor n} \subset A_{nm^2,n}.$$
Since $\phi_{K'}(A_{nm^2,n}) \subseteq J_n$, we can
conclude that $b = s_*(K) \in J_n$ for every
$K \in {\cal K}_n$. which is equivalent to the desired result
that $s_L(J_n) \subseteq J_n$, which proves Theorem~\ref{thvass}.

\section{Quantum Lie group invariants}

In this section, we argue Theorem~\ref{thquantum}. Let $\frak g$ be a simple
Lie group.  Following the well-known theory of quantum topological invariants
\cite{Turaev}, $U_q({\frak g})$ is a quasitriangular Hopf algebra, and there
corresponds a polynomial invariant of oriented, framed graphs in $S^3$ whose
edges are labelled by finite-dimensional representations of $U_q({\frak g})$
and whose vertices are labelled by invariant tensors over the representations
of the incident edges. In this theory, parallel edges labelled by
representations $V_1,\ldots,V_n$ are equivalent to one edge labelled by $V_1
\tensor \ldots \tensor V_n$, and reversing the orientation of an edge is
equivalent to taking the dual of its representation.  The invariants are also
multiadditive under direct sums.

The standard theory can be generalized in a simple way for links so that each
link component is labelled by a formal trace $\tau$ on $U_q({\frak g})$ such
that $\tau(ab) = \tau(ba)$.  If $\tau$ is the character of a linear
representation, then the generalized invariant is the same as before.  The
theory is then multilinear over formal traces, which form a vector space
$T_q({\frak g})$. In particular, given a knot $K$, there is a universal
invariant $Q_{q,{\frak g}}(K) \in T_q({\frak g})^*$.  Moreover,
the characters of irreducible representations form a basis for
$T_q({\frak g})$ since $U_q(\frak g)$ is semisimple for $q$ an
indeterminate, so the generalization is in this case nothing
more than the linear completion of the invariants in which edges
are labelled by true representations.

Given a 2-component link $L$ as in the previous section, if we interpret the
second component $L_2$ as a knot in a solid torus and label it by a
representation $V$, then in the quantum topological invariant theory it
represents some equivariant endomorphism of $V_1 \tensor \ldots \tensor V_n$,
where each $V_i$ is either $V$ or $V^*$ (see Figure~\ref{fendo}).
The tensor product $V = V_1 \tensor \ldots \tensor V_i$ has a direct sum
decomposition
$$V = \bigoplus_{W \in {\cal R}} m_W W,$$
where $\cal R$ a set with one representative of each irreducible
representation and $m_W$ is a non-negative integer.  The map
$L_2$, as an endomorphism of $V$, decomposes as
$$L_2 = \bigoplus L_W \tensor I_W,$$
where $L_W$ is an $m_W \times m_W$ matrix over the ground field and
$I_W$ is the identity on $W$.  If $\tau_{q,\frak g}(K,V)$
is the quantum invariant of $K$ colored by the representation $V$,
then
\begin{equation}
\tau_{q,\frak g}(s_L(K),V) = \sum \Tr(s_W)\tau_{q,\frak g}(K,W) \label{eqsat}
\end{equation}
The left side is a satellited quantum group invariant of the knot $K$.
The right side is a linear combination of ordinary quantum group
invariants with different representations.  Equation~\ref{eqsat}
establishes that $\tau_{q,\frak g}(s_L(K),V)$ factors
through $Q_{q,{\frak g}}(K)$.  Therefore $s_L$ intertwines
$Q_{q,\frak g}$.

\fig{}{A tangle inducing a linear endomorphism.}

That quantum group invariants are invariant under inversion is
equivalent to the statement that
$$\tau_{q,\frak g}(K,W) = \tau_{q,\frak g}(K,W^*)$$
for every $W$.  Recall that the quantum group $U_q({\frak g})$ is constructed
from the Dynkin diagram of $\frak g$, and that any such Dynkin diagram has a
dualizing automorphism, \ie, an automorphism which takes the highest weight
of an irreducible representation $W$ to the highest weight of the dual
representation $W^*$.  It follows that the pair $(U_q({\frak g}),W)$ is
isomorphic to the pair $(U_q({\frak g}),W^*)$, and that the invariants
are equal.

\section{Finite group invariants}

The treatment of quantum group invariants generalizes from $U_q({\frak g})$ to
arbitrary quasitriangular Hopf algebras $H$. However, in the general case,
there is no guarantee that a dualizing homomorphism exists. In particular, if
$H = D(\C[G])$ is the quantum double of a finite group algebra, the
corresponding invariant $\tau_H$ counts group homomorphisms from $\pi_1(S^3 -
K)$ to $G$ with a specified restriction to the peripheral subgroup.  Strictly
speaking, one fixes a group element $g$ and a linear representation $V$ of
the centralizer of $g$, and one takes the sum of $\Tr_V(l)$ over
homomorphisms that take $m$ to $g$, where $m$ is the meridian and $l$ is the
longitude.  However, after a transformation, this is the same as counting
homomorphisms that take $(m,l)$ to $(g,h)$ for fixed $g$ and $h$.  Although
it is interesting to consider this invariant in quantum terms
\cite{Dijkgraaf-Witten,Kuperberg:Invol}, a more sober point of view in
practice is that of classical group theory and algebraic topology.

Define a group-theoretic knot to be a triple $(G,m,l)$, where $G$ is a group
and $m$ and $l$ are two arbitrary commuting elements of $G$ called the
meridian and the longitude.  An honest knot $K$ in $S^3$ yields a
group-theoretic knot $(\pi_1(S^3-K),m,l)$ which determines $K$ completely. 
Define $(G,m,l)$ to be {\em invertible} if there exists an automorphism of
$G$ that inverts $m$ and $l$. This too is equivalent to invertibility for
honest knots. If the number of homomorphisms from $(\pi_1(S^3 - K),m_K,l_K)$
to $(G,m_G,l_G)$ differs from the number of homomorphisms to
$(G,m_G^{-1},l_G^{-1})$, then $K$ is clearly non-invertible.  Moreover, by
inclusion-exclusion over subgroups of $G$, counting homomorphisms yields
the same information as counting epimorphisms.

Let $C_n$ be the $n$th commutator subgroup of $\pi = \pi_1(S^3 - K)$. The
first interesting group-theoretic model of $K$ is $\pi/C_2$, whose structure
is characterized by the Alexander polynomial $\Delta_K(t)$. Unfortunately,
the well-known identity $\Delta_K(t) = \Delta_K(t^{-1})$ implies that
$\pi/C_2$ is an invertible group-theoretic knot.  In other words, all
3-dimensional knots are invertible at the metabelian level; if $G$ is
metabelian, counting homomorphisms to $G$ cannot detect invertibility.

The quotient $\pi/C_3$ is more fruitful:  Hartley \cite{Hartley} found that
the homology of metabelian coverings often does detect invertibility, which
amounts to counting homomorphisms to certain meta-metabelian groups (groups
with metabelian commutator). However, neither these groups nor any other
solvable groups will yield any information when the Alexander polynomial of
$K$ is trivial. In this case, $C_1$ is a perfect group, which means that
$\pi/C_n \cong \Z$ for all $n$, which is the same as the unknot, which is
invertible. The first example of such a knot is Conway's 11-crossing knot,
shown in Figure~\ref{fconway}.

\fig{}{Conway's knot, which has trivial Alexander polynomial.}

The sporadic simple group $M_{11}$ of order 7920 admits no automorphism that
inverts an element $g$ of order 11.  We report here that, if $K$ is the
Conway knot, then $\pi_1(S^3-K)$ admits precisely one epimorphism to $M_{11}$
that sends the meridian $m$ to either $g$ or $g^{-1}$. It follows that $K$ is
non-invertible.  This result was obtained by a computer search:  working with
the Wirtinger presentation of the knot group, the algorithm finds among all
maps from the generating set to the conjugacy class of $g$ those which
satisfy the relations.  Note that Wirtinger generators are conjugate because
they all meridians, so we need not consider other maps from the generators to
$M_{11}$. The first generator can be mapped to $g$ itself, and if the
generators are properly ordered, the first three determine the others in the
case of the Conway knot.  Those who wish to investigate the validity of this
computational result can find the author's software on the Internet
\cite{Kuperberg:m11}.

SnapPea can also prove that $K$ is non-invertible using the
hyperbolic structure of $S^3 - K$ and Mostow rigidity,
but the proof using $M_{11}$ is perhaps more elementary, although
it does involve computer calculations.

Kenichi Kawagoe \cite{Kawagoe:personal} has verified the author's
computational results, and has also verified by computer that the
finite groups $\Sz(8)$ and $\Aut(U_3(3))$ can be used to determine the
non-invertibility of both the Conway knot and the Kinoshita-Terasaka
knot, which is a mutant of the Conway knot.  These results support the
hypothesis that a knot group, like a free group, has a good chance
of admitting any particular finite group $G$ as a quotient, provided
that $G$ is generated by sufficiently few elements. While such a
quotient may prove useful as an invariant, it does not necessarily
reveal any special connection between the knot and the group $G$.


\end{document}